\documentclass[a4paper,twoside,fleqn,notitlepage]{article}%
\usepackage{amsmath}
\usepackage{amsfonts}
\usepackage{amssymb}
\usepackage{graphicx}%
\providecommand{\U}[1]{\protect \rule{.1in}{.1in}}
\newtheorem{theorem}{Theorem}

\newtheorem{corollary}[theorem]{Corollary}

\newtheorem{definition}[theorem]{Definition}

\newtheorem{lemma}[theorem]{Lemma}

\newtheorem{proposition}[theorem]{Proposition}

\newenvironment{proof}[1][Proof]{\noindent \textbf{#1.} }{\  \rule{0.5em}{0.5em}}
\begin{document}

\title{Study of New Class of $q$-fractional Derivative Operator and its Properties}
\author{M. Momenzadeh, S. Norouzpoor\\Near East University\\Lefkosa, TRNC, Mersin 10, Turkey \\Istanbul Gelisim University, Istanbul, Turkey\\Email: mohammad.momenzadeh@neu.edu.tr\\snorouzpoor@gelisim.edu.tr\\ \  \  \  \  \  \  \  \  \ }
\date{}
\maketitle

\begin{abstract}
There are several approaches to the fractional differential operators.
Generalized $q$-fractional integral operator was defined in the aid of
$q$-iterated Cauchy integral and $q$-calculus techniques. In the following
paper, we introduce the Caputo type derivative related to the mentioned
operator and some properties of this operator is investigated. Cauchy problem
based on this operator is studied and existence and uniqueness is discussed.

\end{abstract}

\section{Introduction}

Since 1695, that a letter related to fractional derivative was written,
fractional calculus has created. There are several types of fractional
integral and derivative operators that arise from different aspects. Recent
appplications of fractional differential equations, in explaining natural
phenomena, motivate more and more scientifics to work in this area. One
approach to fractional integral operator is using Cauchy integral. Authors in
\cite{mami} have used $q$-calculus techniques to develop $q$-fractional
difference equations. In this paper, we first state some definitions and
concepts of $q$-calculus, then, we discuss about the concave operator on
Banach space and then, the fixed point theorem will be written based on these
conditions. Some properties of $q$-fractional derivative, integral and Caputo
type of derivative are discussed.

In the first section, let us introduce some familar concepts of q-calculus
that are in use in this paper \cite{kac} and \cite{ernst}. We use $[n]_{q}$ as
a $q-$analogue of any complex number. Naturally, we can define $\left[
n\right]  _{q}!$ as
\[
\left[  a\right]  _{q}=\frac{1-q^{a}}{1-q}\  \  \  \left(  q\neq1\right)
;\  \  \  \left[  0\right]  _{q}!=1;\  \  \  \  \left[  n\right]  _{q}!=\left[
n\right]  _{q}\left[  n-1\right]  _{q}\  \  \  \ n\in \mathbb{N},\  \ a\in
\mathbb{C}\text{ }.
\]
The $q$-shifted factorial and $q$-polynomial coefficient are defined by
\begin{align*}
\left(  a;q\right)  _{0}  &  =1,\  \  \  \left(  a;q\right)  _{n}=%
{\displaystyle \prod \limits_{j=0}^{n-1}}
\left(  1-q^{j}a\right)  ,\  \  \ n\in \mathbb{N},\\
\left(  a;q\right)  _{\infty}  &  =%
{\displaystyle \prod \limits_{j=0}^{\infty}}
\left(  1-q^{j}a\right)  ,\  \  \  \  \left \vert q\right \vert <1,\  \ a\in
\mathbb{C}.
\end{align*}%
\[
\left(
\begin{array}
[c]{c}%
n\\
k
\end{array}
\right)  _{q}=\frac{\left(  q;q\right)  _{n}}{\left(  q;q\right)
_{n-k}\left(  q;q\right)  _{k}},
\]
Let the function $|f(x)x^{\alpha}|$ be bounded on the interval $\left(
0,A\right]  $, for some $0\leq \alpha<1,$ then Jakson integral is defined as
\cite{kac}%

\[
\int \mathit{f}(x)d_{q}x=(1-q)x\sum_{i=0}^{\infty}q^{i}f(q^{i}x)
\]
and it converges to a function $F(x)$ on $\left(  0,A\right]  ,$ which is a
$q-$antiderivative of $f(x)$. Suppose $0<a<b$, the mentioned $q-$integral is
defined as%

\begin{align*}
\int \limits_{0}^{b}\mathit{f}(x)d_{q}x  &  =(1-q)b\sum_{i=0}^{\infty}%
q^{i}f(q^{i}b)\\
\int \limits_{a}^{b}\mathit{f}(x)d_{q}x  &  =\int \limits_{0}^{b}\mathit{f}%
(x)d_{q}x-\int \limits_{0}^{a}\mathit{f}(x)d_{q}x
\end{align*}

\bigskip In addition, we can interchange the order of double $q$-integral by%

\[%
{\displaystyle \int \limits_{0}^{x}}
{\displaystyle \int \limits_{0}^{v}}
d_{q}sd_{q}v=%
{\displaystyle \int \limits_{0}^{x}}
{\displaystyle \int \limits_{qs}^{x}}
d_{q}vd_{q}s
\]

\ Let $C_{q}^{n}[a,b]$ denotes the space of all continues functions with
continuous q-derivatives up to order $n-1$ on the interval $[a,b]$ .
Associated norm function of $C_{q}^{n}[a,b]$ is defined as \cite{mansour}%

\[
\left \Vert f\right \Vert =\sum_{i=0}^{n-1}\max_{a\leq x\leq b}\left \vert
(D_{q}^{i}f)(x)\right \vert \text{ \ ,}f\in C_{q}^{n}[a,b]
\]

The generalize $q$-exponent expression is defined \cite{mami}%

\begin{equation}
(x-y)_{q^{p}}^{\left(  \alpha \right)  }=\prod_{k=0}^{\infty}\frac{\left(
x-y\left(  q^{p}\right)  ^{k}\right)  }{\left(  x-y\left(  q^{p}\right)
^{k+\alpha}\right)  }=\frac{x^{\alpha}\left(  \frac{y}{x};q^{p}\right)
_{\infty}}{\left(  q^{\alpha \left(  p\right)  }\frac{y}{x};q^{p}\right)
_{\infty}}%
\end{equation}
where the normal definition can be expressed as \cite{Zhang}\cite{rajkovic}%

\begin{equation}
(x-a)^{\left(  \alpha \right)  }=x^{\alpha}\prod_{k=0}^{\infty}\frac{\left(
1-\frac{x}{a}q^{k}\right)  }{\left(  1-\frac{x}{a}q^{k+\alpha}\right)  }%
=\frac{\left(  \frac{x}{a};q\right)  _{\infty}}{\left(  q^{\alpha}\frac{x}%
{a};q\right)  _{\infty}}%
\end{equation}
The $q$-derivative of this expression with respect to $x$ and $y$ can be
written as%

\begin{align}
_{x}D_{q}\left[  (x^{p}-y^{p})_{q^{p}}^{\left(  \alpha \right)  }\right]   &
=x^{p-1}\left[  p\alpha \right]  _{q}(x^{p}-y^{p})_{q^{p}}^{\left(
\alpha-1\right)  }\\
_{y}D_{q}\left[  (x^{p}-y^{p})_{q^{p}}^{\left(  \alpha \right)  }\right]   &
=-y^{p-1}\left[  p\alpha \right]  _{q}(x^{p}-\left(  yq\right)  ^{p})_{q^{p}%
}^{\left(  \alpha-1\right)  }%
\end{align}
respectively.

One can express the $q$-Gamma function by using this definition as\cite{gamma
kac}
\begin{equation}
\Gamma_{q}(t)=\frac{(1-q)^{\left(  t-1\right)  }}{(1-q)^{t-1}}%
\end{equation}
\newpage Let us bring the useful lemma proved in \cite{mami}, which express
following integral for $\alpha$ and $\lambda>-1$ ,%

\begin{equation}
\int \limits_{a}^{x}t^{p-1}(x^{p}-\left(  qt\right)  ^{p})_{q^{p}}^{\left(
\alpha-1\right)  }(t^{p}-a^{p})_{q^{p}}^{\left(  \lambda \right)  }d_{q}%
t=\frac{1}{\left[  p\right]  _{q}}\left(  \frac{\Gamma_{q^{p}}(\alpha
)\Gamma_{q^{p}}(\lambda+1)}{\Gamma_{q^{p}}(\alpha+\lambda+1)}\right)  \left[
(x^{p}-a^{p})_{q^{p}}^{\left(  \alpha+\lambda \right)  }\right]
\end{equation}
There are several approaches to fractional differential operators. One of the
demonstration methods of fractional differential equation is using the
itterated Cauchy integrals. In \cite{mami}, authors calculated%

\[
\int \limits_{0}^{x}\left(  \mathit{t}_{1}\right)  ^{p-1}\mathit{d}%
_{q}\mathit{t}_{1}\int \limits_{a}^{\mathit{t}_{1}}\left(  \mathit{t}%
_{1}\right)  ^{p-1}\mathit{d}_{q}\mathit{t}_{2}...\int \limits_{a}%
^{\mathit{t}_{n-1}}\left(  \mathit{t}_{n}\right)  ^{p-1}f(t_{n})\mathit{dt}%
_{n}=\frac{1}{\prod_{k=1}^{n-1}\left[  kp\right]  _{q}}\int \limits_{0}%
^{x}w^{p-1}f(w)\prod_{k=0}^{n-1}\left(  x^{p}-\left(  wq\right)  ^{p}%
q^{kp}\right)  \mathit{d}_{q}\mathit{w}%
\]
In the aid of this calculation, generalized $q$-fractional difference integral
operator is defined by%

\begin{align}
J_{p,q}^{\alpha}\left(  f(x)\right)   &  =\frac{\left(  1-q\right)
^{\alpha-1}}{\left(  1-q^{p}\right)  _{q^{p}}^{\left(  \alpha-1\right)  }}%
\int \limits_{0}^{x}w^{p-1}f(w)(x^{p}-\left(  wq\right)  ^{p})_{q^{p}}^{\left(
\alpha-1\right)  }\mathit{d}_{q}\mathit{w}\\
&  \mathit{=}\frac{\left(  \left[  p\right]  _{q}\right)  ^{1-\alpha}}%
{\Gamma_{q^{p}}(\alpha)}\int \limits_{0}^{x}w^{p-1}f(w)(x^{p}-\left(
wq\right)  ^{p})_{q^{p}}^{\left(  \alpha-1\right)  }\mathit{d}_{q}\mathit{w}%
\end{align}
The corresponding inverse operator named as $q$-fractional difference operator
is defined by%

\begin{align}
(D_{p,q}^{0}f)(x)  &  =f(x)\\
(D_{p,q}^{\alpha}f)(x)  &  =(x^{1-p}D_{q})^{n}\left(  J_{p,q}^{n-\alpha
}\right)  f(x)=\frac{\left(  \left[  p\right]  _{q}\right)  ^{\alpha-n+1}%
}{\Gamma_{q^{p}}(n-\alpha)}(x^{1-p}D_{q})^{n}\int \limits_{0}^{x}%
w^{p-1}f(w)(x^{p}-\left(  wq\right)  ^{p})_{q^{p}}^{\left(  n-\alpha-1\right)
}\mathit{d}_{q}\mathit{w}\text{ \ }%
\end{align}
for $\alpha \geq0$ and $n=\left \lfloor \alpha \right \rfloor +1$ and $p>0$.

\begin{definition}
Let $0<a<b<\infty$, $f:[a,b]\rightarrow R$ be a $q$-integrable function, and
let $\alpha \in(0,1)$ and \ $p>0$ be two fixed real numbers. The Caputo type
$q-$ fractional derivative of order $\alpha$ is defined by%
\begin{equation}
(^{c}D_{a^{+},p,q}^{\alpha}f)(t)=(D_{a^{+},p,q}^{\alpha})(f(t)-f(a))=\frac
{\left(  \left[  p\right]  _{q}\right)  ^{\alpha}}{\Gamma_{q^{p}}(1-\alpha
)}(x^{1-p}D_{q})\int \limits_{a}^{x}w^{p-1}\left[  f(w)-f(a)\right]
(x^{p}-\left(  wq\right)  ^{p})_{q^{p}}^{\left(  -\alpha \right)  }%
\mathit{d}_{q}\mathit{w}\text{ }%
\end{equation}

\end{definition}

Note that, when $f(a)=0,$ the Caputo and the $q-$ fractional derivatives
coincide. Moreover, there exists a relation between these two types of
$q-$fractional derivatives. Namely, if both types of $q-$fractional derivative
exist, then by considering (3) and (4) we can see that,%

\begin{equation}
(^{c}D_{a^{+},p,q}^{\alpha}f)(x)=(D_{a^{+},p,q}^{\alpha})(f(x))-(D_{a^{+}%
,p,q}^{\alpha})(f(a))=(D_{a^{+},p,q}^{\alpha})(f(x))-\frac{\left[
f(0)\right]  \left(  \left[  p\right]  _{q}\right)  ^{\alpha}}{\Gamma_{q^{p}%
}(1-\alpha)}(x^{p}-a^{p})_{q^{p}}^{\left(  -\alpha \right)  }\text{ }%
\end{equation}

Now, consider the simplification os the Caputo $q$-derivative in the next proposition;

\begin{proposition}
The Caputo $q$-derivative can be simplified and written as,
\begin{equation}
(^{c}D_{a^{+},p,q}^{\alpha}f)(x)=\frac{\left(  \left[  p\right]  _{q}\right)
^{\alpha}}{\Gamma_{q^{p}}(1-\alpha)}\int \limits_{a}^{x}\left[  D_{q}%
f(w)\right]  (x^{p}-w^{p})_{q^{p}}^{\left(  -\alpha \right)  }\mathit{d}%
_{q}\mathit{w}\text{ }%
\end{equation}
if all $q$-Jackson integrals in (10) be convergent.

\begin{proof}
By using definition (1) and identification (4) we can get;
\[
(^{c}D_{a^{+},p,q}^{\alpha}f)(x)=\frac{\left(  \left[  p\right]  _{q}\right)
^{\alpha}}{\left[  p(1-\alpha)\right]  _{q}\Gamma_{q^{p}}(1-\alpha)}%
(x^{1-p}D_{q})\int \limits_{a}^{x}\left(  f(w)-f(a)\right)  \left[  D_{q}%
(x^{p}-w^{p})_{q^{p}}^{\left(  1-\alpha \right)  }\right]  \mathit{d}%
_{q}\mathit{w}%
\]

The function $(x^{p}-w^{p})_{q^{p}}^{\left(  1-\alpha \right)  }$ is q-regular
and we can apply q-integral by part to reach
\[
(^{c}D_{a^{+},p,q}^{\alpha}f)(x)=\frac{\left(  \left[  p\right]  _{q}\right)
^{\alpha}}{\left[  p(1-\alpha)\right]  _{q}\Gamma_{q^{p}}(1-\alpha)}%
(x^{1-p}D_{q})\int \limits_{a}^{x}(x^{p}-\left(  wq\right)  ^{p})_{q^{p}%
}^{\left(  1-\alpha \right)  }\left[  D_{q}f(w)\right]  \mathit{d}%
_{q}\mathit{w}%
\]

In the aid of lemma 1.12 of \cite{mansour}, we can apply the differentiation
inside of Jackson integral. In fact, the derivative respect to $x$ is defined
at (3) and $(\left(  xq\right)  ^{p}-\left(  xq\right)  ^{p})_{q^{p}}^{\left(
1-\alpha \right)  }=0$ which satisfied the assumption of lemma 1.12, so
applying (3) will complete the proof.

\begin{corollary}
The direct consequence of this proposition can demonstrate the relation
between Caputo type and normal fractional derivative operators (10) and (11)
as follow%

\begin{align}
(^{c}D_{a^{+},p,q}^{\alpha}f)(x)  &  \mathit{=}\frac{1}{\left[  p(1-\alpha
)\right]  _{q}}(x^{1-p}D_{q})\left(  J_{a^{+},p,q}^{1-\alpha}\left(
w^{1-p}D_{q}\left(  f\right)  \right)  \right)  (x)\\
&  \mathit{=}\frac{1}{\left[  p(1-\alpha)\right]  _{q}}(D_{a^{+},p,q}^{\alpha
})\left(  w^{1-p}D_{q}\left(  f\right)  \right)  (x)
\end{align}

\end{corollary}
\end{proof}
\end{proposition}

In fact, study of fractional differential operator for the case $\alpha
\in \left(  0,1\right)  $ is more attractive and serves many purposes. First,
Abel integral equation which is inspired by a problem of mechanics is one of
the first integral equations ever treated and considered the fractional
integral operator for $\alpha$ between zero and one.\cite{abel} As a second
main reason, we can mention the application of fractional differential
operator for modelling different phenomena that were used $\alpha \in \left(
0,1\right)  $ commonly to describe simulation. \cite{cu} In addition, the
calculations for the case $\alpha \in \left(  0,1\right)  $ is not too long and
complicated and also most of them can be generalized for arbitrary $\alpha.$
Our porpose in this article is introducing this operator clearly such that can
be used later for application cases.

On the other hand, Kilbas et. al. substitute $f(a)+%
{\displaystyle \int \limits_{a}^{t}}
f^{\prime}(w)dw$ instead of $f(t)$ to find the formula for derivative operator
and this technique can be used for (11) as well to reach $q$-analogue of
it.\cite{kilbas} We can apply foundemental theorem of $q$-calculus to see;%

\begin{align}
(^{c}D_{a^{+},p,q}^{\alpha}f)(x)  &  =(D_{a^{+},p,q}^{\alpha}%
)(f(x)-f(a))=\frac{\left(  \left[  p\right]  _{q}\right)  ^{\alpha}}%
{\Gamma_{q^{p}}(1-\alpha)}(x^{1-p}D_{q})\int \limits_{a}^{x}w^{p-1}\left[
{\displaystyle \int \limits_{a}^{w}}
D_{q}f(t)d_{q}t\right]  (x^{p}-\left(  wq\right)  ^{p})_{q^{p}}^{\left(
-\alpha \right)  }\mathit{d}_{q}\mathit{w}\text{ }\\
&  =\frac{\left(  \left[  p\right]  _{q}\right)  ^{\alpha}}{\Gamma_{q^{p}%
}(1-\alpha)}(x^{1-p}D_{q})\int \limits_{a}^{x}%
{\displaystyle \int \limits_{qt}^{x}}
\left[  D_{q}f(t)\right]  w^{p-1}(x^{p}-\left(  wq\right)  ^{p})_{q^{p}%
}^{\left(  -\alpha \right)  }\mathit{d}_{q}\mathit{w}d_{q}t\text{ }\\
&  =\frac{\left(  \left[  p\right]  _{q}\right)  ^{\alpha}}{\left[
p(1-\alpha)\right]  _{q}\Gamma_{q^{p}}(1-\alpha)}(x^{1-p}D_{q})\int
\limits_{a}^{x}(x^{p}-\left(  qt\right)  ^{p})_{q^{p}}^{\left(  1-\alpha
\right)  }\left[  D_{q}f(t)\right]  \mathit{d}_{q}\mathit{t}%
\end{align}

In (17), we interchange the order of iterated $q$-integral and reach to
similar result. We can continue this procedure to find more general identity
which leads to characterization of the classes of function which fractional
operator (7) is defined for them.

\section{Some Properties of Integral and Derivative Operator}

In this section, some properties of Caputo-type derivative and integral
operator are studied. These propeties lead to solving related $q$-fractional
difference equation. First, we should state the boundness of $q-$fractional
integral operator.

\begin{proposition}
Let $f(x)$ be a real function which is defined on $q$-geometric set or
particularly assume that $x$ can be formed as $bq^{n}$, then $q$-fractional
integral operator $J_{a^{+},p,q}^{\alpha}$ is linear and bounded from
$C_{q}[a,b]$ to $C_{q}[a,b]$, that means%
\begin{equation}
\left \Vert J_{a^{+},p,q}^{\alpha}\left(  f\right)  (x)\right \Vert \leq
A_{p,q,\alpha}\left \Vert f(x)\right \Vert
\end{equation}

\begin{proof}
Assume that $f\in C_{q}[a,b]$ that the space is equipped with maximum norm. In
addition, assume that $x$ can be formed as $bq^{n}$ for some $n\in%
\mathbb{N}
,$ then we can apply $q$-analogue of H\"{o}lder inequality as\cite{mansour}%

\begin{align*}
\left \vert J_{a^{+},p,q}^{\alpha}\left(  f\right)  (x)\right \vert  &
\leq \frac{\left(  \left[  p\right]  _{q}\right)  ^{1-\alpha}}{\Gamma_{q^{p}%
}(\alpha)}\left \Vert f\right \Vert _{C}\left \vert \int \limits_{a}^{x}%
w^{p-1}(x^{p}-\left(  wq\right)  ^{p})_{q^{p}}^{\left(  \alpha-1\right)
}\mathit{d}_{q}\mathit{w}\right \vert \\
&  =\frac{\left(  \left[  p\right]  _{q}\right)  ^{1-\alpha}}{\left[
p\alpha \right]  _{q}\Gamma_{q^{p}}(\alpha)}\left \Vert f\right \Vert
_{C}\left \vert (x^{p}-a^{p})_{q^{p}}^{\left(  \alpha \right)  }\right \vert
\end{align*}

Here, we used identity (4) to calculate the $q$-integral. Since Jackson
integral is a linear operator, linearlity is obvious.
\end{proof}
\end{proposition}

Let us consider the inverse operator of integral operator, that we introduced
it before. In the next proposition, we show that how the Caputo type
derivative can be operate on integral operator. Semi-group property of
integral operator is investigated\cite{mami}.

\begin{proposition}
Let $0<a<b<\infty$ and $f$ is defined in all following expressions, also
assume that $\alpha \in(0,1)$ and \ $p>0$ be two fixed reals such that
derivative and integral operator are defined, then following identities are true:%

\begin{align}
^{c}D_{a^{+},p,q}^{\alpha}\left(  (J_{a^{+},p,q}^{\alpha}f)(x)\right)   &
=f(x)\\
J_{a^{+},p,q}^{\alpha}\left(  (^{c}D_{a^{+},p,q}^{\alpha}f)(x)\right)   &
=f(x)-f(a)
\end{align}

\begin{proof}
proof of the first part is based on interchanging of integral and applying the
definitions. This is proved at \cite{mami} for derivative opertaor and the
procedure is similar for this operator. We prove the second identity by using
corollary (3) and semi group property of integral operator, so%

\[
J_{a^{+},p,q}^{\alpha}\left(  (^{c}D_{a^{+},p,q}^{\alpha}f)(x)\right)
=J_{a^{+},p,q}^{\alpha}J_{a^{+},p,q}^{1-\alpha}\left(  w^{1-p}D_{q}\left(
f\right)  \right)  (x)=J_{a^{+},p,q}^{1}\left(  w^{1-p}D_{q}\left(  f\right)
\right)  (x)=\int \limits_{0}^{x}D_{q}f(w)\mathit{d}_{q}\mathit{w=}f(x)-f(a)
\]

\end{proof}
\end{proposition}

\section{Cauchy problem of general q-fractional operator}

In this section, Cauchy problem related to the introduced Caputo type
derivative is investigated. We construct solution of this $q$-fractional
difference equation in the aid of discussed properties in last section.
Moreover we apply a fixed point theorem based on $q$-successive Approximations
to guarantee existance and uniqueness of this difference equation. We start by
direct solution of $q$-difference equation in following lemma:

\begin{lemma}
\bigskip let $0<\zeta<1$ and $\alpha \in \left(  0,1\right)  $, if $f(t)$ is a
continuous function on $\left[  0,\infty \right)  \times \left(  \zeta
-r,\zeta+r\right)  $ where $r$ is a positive real constant and $^{c}%
D_{a^{+},p,q}^{\alpha}u$ is defined $,$ then solution of the following
boundary value problem
\begin{align}
(^{c}D_{a^{+},p,q}^{\alpha}u)(t)  &  =f(t,u(t))\text{ \  \  \  \  \  \ }%
0<a<t<b<\infty \\
u(a)  &  =\zeta
\end{align}

can be written as
\begin{equation}
u(t)=\zeta+\frac{\left(  \left[  p\right]  _{q}\right)  ^{1-\alpha}}%
{\Gamma_{q^{p}}(\alpha)}\int \limits_{a}^{t}w^{p-1}f(w,u(w))(t^{p}-\left(
wq\right)  ^{p})_{q^{p}}^{\left(  \alpha-1\right)  }\mathit{d}_{q}%
\mathit{w}\text{ }%
\end{equation}

\begin{proof}
Apply $J_{a^{+},p,q}^{\alpha}$ in both sides of (22) and use (21), then use
initial value to reach the solution.
\end{proof}
\end{lemma}

Following theorem shows the conditions on $f(w,u(w))$ to verify existance and
uniqueness for given $q$-fractional difference equation. Moreover,
$q$-sucessive approximation based on solution of (22) is constructed. Main
theorem of this article determines solution and nature of this approximation.
However, the sprit of proof is similar to the discussion on \cite{mansour}. In
the aid of last lemma, we define the sequence of functions as follow
\begin{equation}
\varphi_{n}(t)=%
\genfrac{\{}{.}{0pt}{}{\zeta
\text{\  \  \  \  \  \  \  \  \  \  \  \  \  \  \  \  \  \  \  \  \  \  \  \  \  \  \  \  \  \  \  \  \  \  \  \  \  \  \  \  \  \  \  \  \  \  \  \  \  \  \  \  \  \  \  \  \  \  \  \  \  \  \  \  \  \  \  \  \  \  \  \  \  \  \  \  \  \  \  \  \  \  \  \  \  \  \  \  \  \  \  \  \ for
\ }n=1}{\zeta+\frac{\left(  \left[  p\right]  _{q}\right)  ^{1-\alpha}}%
{\Gamma_{q^{p}}(\alpha)}\int \limits_{a}^{t}w^{p-1}f(w,\varphi_{n-1}%
(w))(t^{p}-\left(  wq\right)  ^{p})_{q^{p}}^{\left(  \alpha-1\right)
}\mathit{d}_{q}\mathit{w}\text{\  \  \  \  \  \  \  \  \  \ for \ }n\geq2}%
\text{ }%
\end{equation}

\begin{theorem}
\bigskip Assume that all conditions of last lemma holds true. In addition,
assume that $f(w,\varphi(w))$ be a function which is continuous at $w=a$ and
verify the Lipschitz' consition, i.e.%

\[
\left \vert f(w,\widetilde{y})-f(w,y)\right \vert \leq A\left \vert \widetilde
{y}-y\right \vert
\]
Where $A$ is a positive constant and $0\leq w$ and $\left \vert y-\zeta
\right \vert <r$. Then q-fractional differential equation problem (22) with
given initial values (23), has a unique solution.

\begin{proof}
Applying Lipschitz' condition to see
\begin{align*}
\left \vert f(w,y)\right \vert  &  \leq \left \vert f(a,\zeta)\right \vert
+\left \vert f(w,y)-f(w,\zeta)\right \vert +\left \vert f(w,\zeta)-f(a,\zeta
)\right \vert \\
&  \leq \left \vert f(a,\zeta)\right \vert +rA+1
\end{align*}
We mention that, according to our assumptions $f(w,\varphi(w))$ is continuous
at $(a,\varphi(a))=(a,\zeta)$ and for sufficient small $w<\gamma$, $\left \vert
f(w,\zeta)-f(a,\zeta)\right \vert <1$. Define the nonzero constants $K,h$ to be%
\[
K=:Sup_{0<w<\gamma,\left \vert y-\zeta \right \vert <r}\left \vert
f(w,y)\right \vert
\]
We establish the existance of solution $\left \{  \varphi_{n}(w)\right \}  ,$
using the method of succesive approximations. We consider the sequence that is
defined by (25), these functions are continuous, because
\begin{align*}
\left \vert \varphi_{n}(t)-\zeta \right \vert  &  =\frac{\left(  \left[
p\right]  _{q}\right)  ^{1-\alpha}}{\Gamma_{q^{p}}(\alpha)}\left \vert
\int \limits_{a}^{t}w^{p-1}f(w,\varphi_{n-1}(w))(t^{p}-\left(  wq\right)
^{p})_{q^{p}}^{\left(  \alpha-1\right)  }\mathit{d}_{q}\mathit{w}%
\text{\ }\right \vert \text{\  \  \  \ }\\
&  \leq K\frac{\left(  \left[  p\right]  _{q}\right)  ^{1-\alpha}}{\left[
p\alpha \right]  _{q}\Gamma_{q^{p}}(\alpha)}(t^{p}-a^{p})_{q^{p}}^{\left(
\alpha \right)  }%
\end{align*}
Thus, each $\left \{  \varphi_{n}(w)\right \}  $ is continuous at $w=a$ and (25)
is well-defined. In the same manner, if we apply the induction, then we can see%

\begin{equation}
\left \vert \varphi_{n+1}(t)-\varphi_{n}(t)\right \vert \leq \left(
\frac{\left(  \left[  p\right]  _{q}\right)  ^{1-\alpha}}{\left[
p\alpha \right]  _{q}\Gamma_{q^{p}}(\alpha)}(t^{p}-a^{p})_{q^{p}}^{\left(
\alpha \right)  }\right)  ^{n}A^{n-1}K\text{\ }%
\end{equation}
Now consider following series that is constructed under the above assumption,
i.e. general term of series verifies (26).%
\[
\varphi_{i}(x)+%
{\displaystyle \sum \limits_{l=1}^{\infty}}
\varphi_{l+1}(x)-\varphi_{l}(x)
\]
The relation (26) supports the uniform convergence of the series by
Weierstrass M-test. Therefore, $\varphi_{n}(t)$ uniformly convergence to a
function, called $\varphi(t)\ $and we aclaim that $\varphi(t)$ is the solution
of (22). Indeed, from Lipschitz' condition we have:%
\[
\left \vert f(w,\varphi_{m}(x))-f(w,\varphi(x))\right \vert \leq A\left \vert
\varphi_{m}(x)-\varphi(x)\right \vert
\]
Since the right side of above inequality approaches uniformly to zero as
$m\rightarrow \infty,$ it follows that
\[
\lim_{m\rightarrow \infty}f(w,\varphi_{m}(x))=f(w,\varphi(x))
\]
Now apply this limit on (25) to verify that $\varphi(x)$ is the solution of (22).
\end{proof}
\end{theorem}

\section{Example}

In this section, we obtain an example of $q$-fractional difference equation
which leads to this $q$-analogue of Mittag-leffler function. Indeed, following
initial problem is considered:%

\begin{align}
(^{c}D_{0^{+},p,q}^{\alpha}u)(t)  &  =u(t)\text{\  \  \  \  \  \ }0<t<b<\infty \\
u(0)  &  =1
\end{align}

We can construct $q$-successive approximation in the aid of (25), then
sequence of functions $\varphi_{n}(x)$ can be defined as following recurrence formula%

\begin{equation}
\varphi_{n}(x)=%
\genfrac{\{}{.}{0pt}{}{1\text{\  \  \  \  \  \  \  \  \  \  \  \  \  \  \  \  \  \  \  \  \  \  \  \  \  \  \  \  \  \  \  \  \  \  \  \  \  \  \  \  \  \  \  \  \  \  \  \  \  \  \  \  \  \  \  \  \  \  \  \  \  \  \  \  \  \  \  \  \  \  \  \  \  \  \  \  \  \  \ for
\ }n=1}{1+\frac{\left(  \left[  p\right]  _{q}\right)  ^{1-\alpha}}%
{\Gamma_{q^{p}}(\alpha)}\int \limits_{0}^{x}w^{p-1}\varphi_{n-1}(w)(x^{p}%
-\left(  wq\right)  ^{p})_{q^{p}}^{\left(  \alpha-1\right)  }\mathit{d}%
_{q}\mathit{w}\text{\  \  \  \ for \ }n\geq2}%
\text{ }%
\end{equation}

We can calculate a few terms of $\varphi_{n}(x)$ by using (4) and (6), specificly%

\begin{align*}
\varphi_{1}(x) &  =1\text{\  \  \  \ }\\
\varphi_{2}(x) &  =1+\frac{\left(  \left[  p\right]  _{q}\right)  ^{1-\alpha}%
}{\Gamma_{q^{p}}(\alpha)}\int \limits_{0}^{x}w^{p-1}(x^{p}-\left(  wq\right)
^{p})_{q^{p}}^{\left(  \alpha-1\right)  }\mathit{d}_{q}\mathit{w=}%
1+\frac{\left(  \left[  p\right]  _{q}\right)  ^{1-\alpha}}{\left[
p\alpha \right]  _{q}\Gamma_{q^{p}}(\alpha)}(x^{p}-0^{p})_{q^{p}}^{\left(
\alpha \right)  }=1+\frac{\left(  \left[  p\right]  _{q}\right)  ^{-\alpha}%
}{\Gamma_{q^{p}}(\alpha+1)}(x^{p}-0^{p})_{q^{p}}^{\left(  \alpha \right)  }\\
\varphi_{3}(x) &  =1+\frac{\left(  \left[  p\right]  _{q}\right)  ^{-\alpha}%
}{\Gamma_{q^{p}}(\alpha+1)}(x^{p}-0^{p})_{q^{p}}^{\left(  \alpha \right)
}+\frac{\left(  \left[  p\right]  _{q}\right)  ^{2-2\alpha}}{\left[
p\alpha \right]  _{q}\left(  \Gamma_{q^{p}}(\alpha)\right)  ^{2}}%
\int \limits_{0}^{x}w^{p-1}(w^{p}-0^{p})_{q^{p}}^{\left(  \alpha \right)
}(t^{p}-\left(  wq\right)  ^{p})_{q^{p}}^{\left(  \alpha-1\right)  }%
\mathit{d}_{q}\mathit{w}\\
&  =1+\frac{\left(  \left[  p\right]  _{q}\right)  ^{-\alpha}}{\Gamma_{q^{p}%
}(\alpha+1)}(x^{p}-0^{p})_{q^{p}}^{\left(  \alpha \right)  }+\frac{\left(
\left[  p\right]  _{q}\right)  ^{2-2\alpha}}{\left[  p\alpha \right]
_{q}\left(  \Gamma_{q^{p}}(\alpha)\right)  ^{2}}\frac{1}{\left[  p\right]
_{q}}\left(  \frac{\Gamma_{q^{p}}(\alpha)\Gamma_{q^{p}}(\alpha+1)}%
{\Gamma_{q^{p}}(2\alpha+1)}\right)  \left[  (x^{p}-0^{p})_{q^{p}}^{\left(
2\alpha \right)  }\right]  \\
&  =1+\frac{\left(  \left[  p\right]  _{q}\right)  ^{-\alpha}}{\Gamma_{q^{p}%
}(\alpha+1)}(x^{p}-0^{p})_{q^{p}}^{\left(  \alpha \right)  }+\frac{\left(
\left[  p\right]  _{q}\right)  ^{-2\alpha}}{\Gamma_{q^{p}}(2\alpha+1)}%
(x^{p}-0^{p})_{q^{p}}^{\left(  2\alpha \right)  }%
\end{align*}

{}It is easy to see by induction that general solution $\varphi_{m}(x)$ can be
presented by%

\begin{equation}
\varphi_{m}(x)=%
{\displaystyle \sum \limits_{n=0}^{m}}
\frac{\left(  \left[  p\right]  _{q}\right)  ^{-n\alpha}}{\Gamma_{q^{p}%
}(n\alpha+1)}(x^{p}-0^{p})_{q^{p}}^{\left(  n\alpha \right)  }\text{ }%
\end{equation}

Definitely, by putting all parameter to suitable values, $\varphi_{n}(x)$
approaches to exponential function and form of solution motivates the
definition of Mittag-Leffler. More details about properties and definition of
classical Mittag-Leffler can be found at \cite{sa}\cite{ka} and can be
compared by $\varphi_{n}(x).$

\end{document}